\documentclass[journal]{IEEEtran} 
\IEEEoverridecommandlockouts

\usepackage{graphicx}
\usepackage{amssymb}
\usepackage{amsmath}
\usepackage{epsfig}
\usepackage{epsf}
\usepackage{subfig}
\usepackage[compress]{cite}

\def\BibTeX{{\rm B\kern-.05em{\sc i\kern-.025em b}\kern-.08em
		T\kern-.1667em\lower.7ex\hbox{E}\kern-.125emX}}

\newtheorem{definition}{Definition}[section]

\newtheorem{lemma}{Lemma}[section]

\begin{document}
	
	
	\title{A Graphical-Based Method for Plotting Local Bifurcation Diagram}
	
	\author{Shahram~Aghaei~and~Abolghasem~Daeichian
		\thanks{S. Aghaei is with Electrical and Computer Engineering Department, Yazd University, Yazd, Iran; aghaei@yazd.ac.ir}
		\thanks{A. Daeichian is with Department of Electrical Engineering, Faculty of Engineering, Arak University, Arak, 38156-8-8349 Iran; a-daeichian@araku.ac.ir, a.daeichian@gmail.com}
		\thanks{Please cite thia paper as: Aghaei, Shahram, and Abolghasem Daeichian. "A Step-by-Step Algorithm for Plotting Local Bifurcation Diagram." arXiv: (2021)}}
	
	\maketitle

\begin{abstract}

{\em Bifurcation diagram is a powerful tool that visually gives information about the behavior of the equilibrium points of a dynamical system respect to the varying parameter. This paper proposes an educational algorithm by which the local bifurcation diagram could be plotted manually and fast in an easy and straightforward way. To the students, this algorithmic method seems to be simpler and more straightforward than mathematical plotting methods in educational and ordinary problems during the learning and studying of courses  related to dynamical systems and bifurcation diagrams. For validation, the algorithm has been applied to several educational examples in the course of dynamical systems.}\\

\end{abstract}

\section{Introduction}
Bifurcation theory is mathematical study of changes in the qualitative or topological structure of a given family, such as the integral curves of a family of vector fields, and the solutions of a family of differential equations. A bifurcation occurs when a small smooth change made to the parameter values (the bifurcation parameters) of a system, causes a sudden qualitative or topological change in its behavior \cite{Blanchard}. Bifurcations occur in both continuous systems (described by ODEs, DDEs or PDEs), and discrete systems (described by maps).
It is useful to divide bifurcations into two principal classes. Local bifurcations, which can be analyzed entirely through changes in the local stability properties of equilibria, periodic orbits or other invariant sets as parameters cross through critical thresholds; and Global bifurcations, which often occur when larger invariant sets of the system collide with each other, or with equilibria of the system. They cannot be detected purely by a stability analysis of the equilibria (fixed points) \cite{Perko}.\newline\indent
Local Bifurcation in continuous systems corresponds to the real part of an eigenvalue of an equilibrium passing through zero. In discrete systems (those described by maps rather than ODEs), this corresponds to a fixed point having a Floquet multiplier with modulus equal to one. In both cases, the equilibrium is non-hyperbolic at the bifurcation point. The topological changes in the phase portrait of the system can be confined to arbitrarily small neighbourhood of the bifurcating fixed points by moving the bifurcation parameter close to the bifurcation point (hence "local") \cite{Kocak}. Examples of local bifurcations include Saddle-node (fold) bifurcation, Transcritical bifurcation, Pitchfork bifurcation, Period-doubling (flip) bifurcation, Hopf bifurcation and Neimark (secondary Hopf) bifurcation.
A bifurcation diagram shows the possible long-term values (equilibria/fixed points or periodic orbits) of a system as a function of a bifurcation parameter in the system \cite{Glendinning}. \\
This paper introduces and develops a method, based on the well known root loci algorithm, originally introduced by Evans \cite{Evans1,Evans2}, to plot local bifurcation diagram. As same as root-luci in control system course, this method could be used as a tool for training the section of bifurcation diagrams in dynamical system courses and also as a useful manual tool for plotting such a diagrams. According to the authors experiences in the course of dynamical systems, this algorithm seems to be simpler and more straightforward than the analytical ways of plotting in educational and ordinary problems. Also, It seems to be an efficient technical algorithm for plotting the bifurcation diagrams in teaching of mathematical courses such as dynamical systems. \newline\indent
Section \ref{Problem statements} presents the problem statement. In section \ref{algorithm}, the algorithm is introduced. Section \ref{Examples} validates the proposed method by some examples. The paper is concluded in section \ref{Conclusion}.

\section{Problem statements}
\label{Problem statements}
Consider the following dynamical system:
\begin{equation} \label{system1}
    x^{(n)} = f(x)+\lambda g(x)
\end{equation}
where $x$ is the state of system and $\lambda \in R$ is a variable parameter. The objective is to develop an algorithm to plot the bifurcation diagram respect to $\lambda$. The equilibrium points of system Eq.\ref{system1} are the roots of the following equation:
\begin{equation} \label{system2}
    f(x)+\lambda g(x)= 0
\end{equation}
Obviously, the roots of the Eq.\ref{system2} are depended on the parameter $\lambda$ except those values that are common roots of $f(x)$ and $g(x)$. Thus Eq.\ref{system2} can be rewritten in the form of Eq.\ref{system3}:
\begin{equation} \label{system3}
    h(x)(f_1(x)+\lambda g_1(x))= 0
\end{equation}
where $f_1(x)$ and $g_1(x)$ are coprime and $h(x)$ is the common factor of $f(x)$ and $g(x)$ and its roots are independent of $\lambda$.
\newtheorem{assump}{Remark}
\begin{assump}\label{remark1}
    The roots of $h(x)$ are considered as constant equilibria of Eq.\ref{system1}. So, the following focuses on the roots of
    \begin{equation} \label{system4}
        f_1(x)+\lambda g_1(x)= 0
    \end{equation}
    as equilibria depended on the$\lambda$.
\end{assump}
\begin{lemma}\label{lemma1}
    if $\lambda=0$ then the roots of $f_1(x)$ are the equilibrium points of Eq.\ref{system1}.\\
    \begin{proof}
    considering Eq.\ref{system4} with $\lambda=0$, obviously, the proof is completed.
    \end{proof}
\end{lemma}

\begin{lemma}\label{lemma2}
    if $\lambda=\pm\infty$ then the roots of $g_1(x)$ are the equilibrium points of Eq.\ref{system1}.\\
    \begin{proof}
    dividing Eq.\ref{system4} by $\lambda=0$ yields  $\frac{f_1 (x)}{\lambda}+g_1(x)=0$. Then with $\lambda=\pm\infty$, obviously, the proof is completed.
    \end{proof}
\end{lemma}
\begin{assump}
    Based on root loci \cite{Evans1, Evans2}, the Eq.\ref{system4} can be rewritten as:
    \begin{equation} \label{system5}
        1+\mu \frac{f_1(x)}{g_1(x)}= 0
    \end{equation}
    in which the greatest degree of $f_1 (x)$ and $g_1 (x)$ have positive coefficients and $\mu = \lambda$ or $\mu = -\lambda$. The roots of $f_1 (x)$ and $g_1 (x)$ are called as poles and zeros, respectively.
\end{assump}

\section{Bifurcation Diagram Plotting Algorithm}
\label{algorithm}

Step 1: Following Eqs.\ref{system1} to \ref{system4} then write Eq.\ref{system5}.\newline\indent
Step 2: Locate the poles and zeros on the vertical axis of bifurcation diagram. Sign the sections of vertical axis that above which there are odd number of poles and zeros as locus for $\mu>0$ and the other sections as locus for $\mu <0$.
    \begin{proof}
    write the Eq.\ref{system5}in the form of poles and zeros as follows:
    \begin{equation} \label{system6}
        \frac{\Pi(x-z_i)}{\Pi(x-p_i)}= -\frac{1}{\mu}
    \end{equation}
    where $z_i$ and $p_i$ are zeros and poles, respectively. Thus $\forall x \in R$ if $\Sigma\angle(x-z_i)-\Sigma\angle(x-p_i)=0 \pm 2k \pi $, then there exist a $\mu>0$ that $x$ lies on roots of Eq.\ref{system6}.
    \begin{eqnarray}
        &&\Sigma\angle(x-z_i)-\Sigma\angle(x-p_i)= \nonumber\\
        &&\Sigma_{\forall z_i>x}\angle (x-z_i)-\Sigma_{\forall p_i>x}\angle (x-p_i)+\nonumber\\&&
        \Sigma_{\forall z_i<x}\angle (x-z_i)-\Sigma_{\forall p_i<x}\angle (x-p_i) =\nonumber\\
        &&\Sigma_{\forall z_i>x}\pi-
        \Sigma_{\forall p_i>x}(-\pi)= 0 \pm 2k \pi \nonumber
    \end{eqnarray}
    Providing $m+n=2k$ in which $m$ and $n$ are the numbers of $z_i>x$ and $p_i>x$, respectively. Repeating the above procedure lead to $m+n=2k+1$ for $\mu<0$.
    \end{proof}

Step 3: Consider a vertical linear asymptote at finite value $\mu=-\lim_{x \to \infty}\frac{f_1(x)}{g_1(x)}$, if there is any.\\
    \begin{proof}
        It's clear from Eq.\ref{system5}.
    \end{proof}

Step 4: Consider horizontal linear asymptotes at finite zeros.\\
    \begin{proof}
    It's result of lemma\ref{lemma2}.
    \end{proof}

Step 5: Sketch the branches of bifurcation diagram, starting from poles, and proceed to asymptotes in each section, considering the sign of $\mu$.\\
    \begin{proof}
    See lemmas\ref{lemma1} and \ref{lemma2}.
    \end{proof}

Step 6: Sketch horizontal lines corresponding to the roots of $h(x)$ in Eq.\ref{system3} as the branches of bifurcation diagram independent of $\lambda$ (or equivalently $\mu$).\\
    \begin{proof}
    See remark \ref{remark1}.
    \end{proof}

Step 7: If $\mu=-\lambda$, then flip the bifurcation diagram horizontally.
\begin{definition}
    A branch is a segment of bifurcation diagram along which there is no change in the sign of slope.
\end{definition}

Step 8: If $f(x)+\lambda g(x)>0$ for $x\gg \max\{{z_i,p_i}\}$ then the upper branch is unstable, next one is stable and so on; otherwise the upper branch is stable, next one is unstable and so on.\\
    \begin{proof}
    The proof is trivial considering the diagram of $f(x)+\lambda g(x)$ versus $x$ for any $\lambda$.
    \end{proof}

\section{Examples}
\label{Examples}

\begin{figure}[htb]
\begin{center}$
\begin{array}{ccc}
\includegraphics[width=2.7cm]{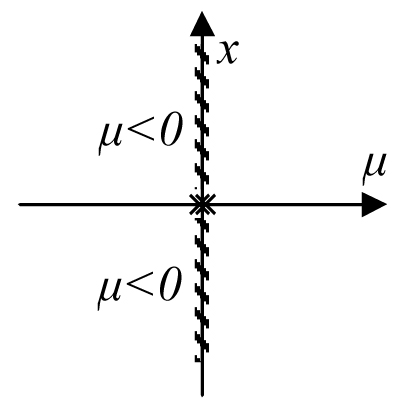}&
\includegraphics[width=2.7cm]{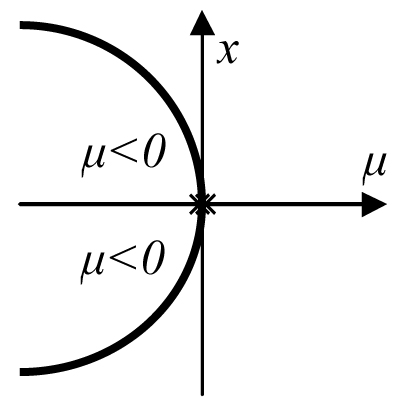}
\includegraphics[width=2.7cm]{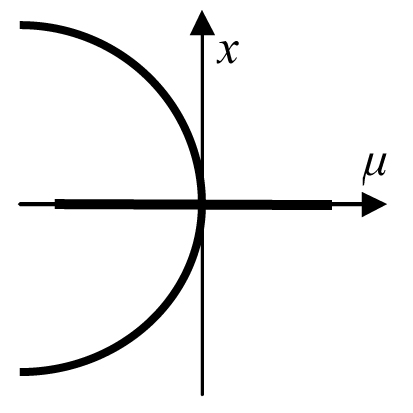}\\
a&b&c
\end{array}$\\
$\begin{array}{cc}
\includegraphics[width=2.7cm]{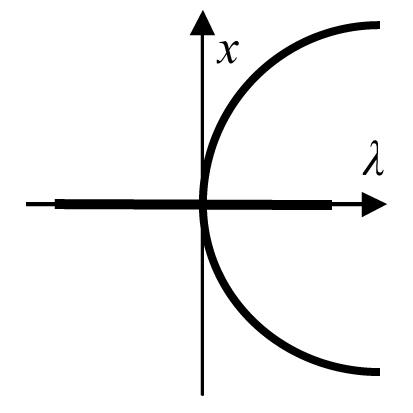}&
\includegraphics[width=2.7cm]{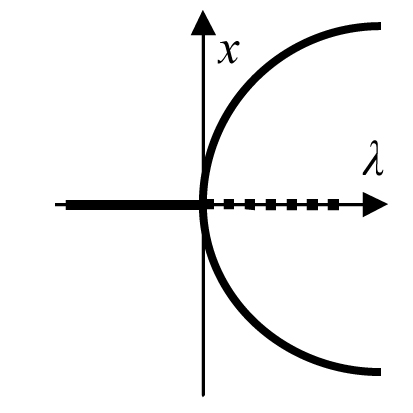}\\
d&e
\end{array}$
\caption{Example 1}
\label{Ex1}
\end{center}
\end{figure}

\newtheorem{Example}{Example}
\begin{Example}
Consider the one-dimensional system $\dot x=\lambda x-x^3$.
\end{Example}
Step1: $x(1-\lambda \frac{1}{x^2})=x(1+\mu \frac{1}{x^2})$, where $\mu=-\lambda$. There are two poles at $x=0$, no zeros, and a constant root at $x=0$.
Step 2: See Fig.\ref{Ex1}.a.
Step 3: There is no finite value vertical asymptote.
Step 4: There is no zeros.
Step 5: See Fig.\ref{Ex1}.b.
Step 6: See Fig.\ref{Ex1}.c.
Step 7:    See Fig.\ref{Ex1}.d.
Step 8: See Fig.\ref{Ex1}.e.

\begin{figure}[h]
\begin{center}$
\begin{array}{ccc}
\includegraphics[width=2.5cm]{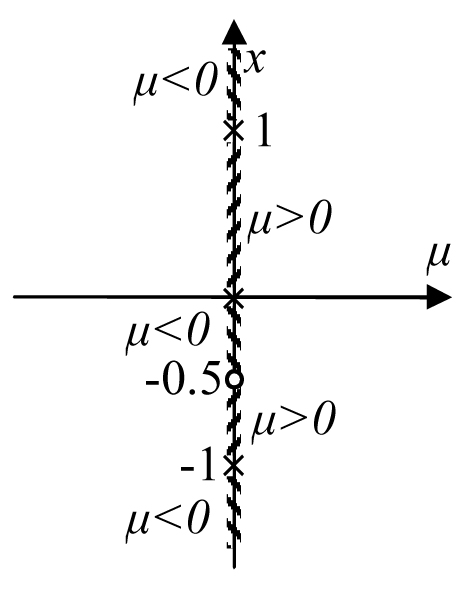}&
\includegraphics[width=2.5cm]{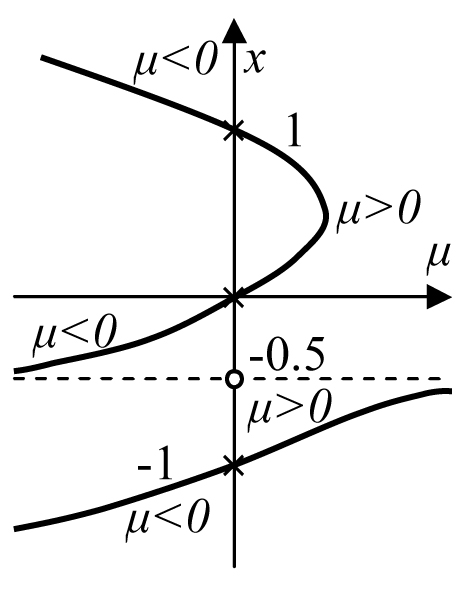}&
\includegraphics[width=2.5cm]{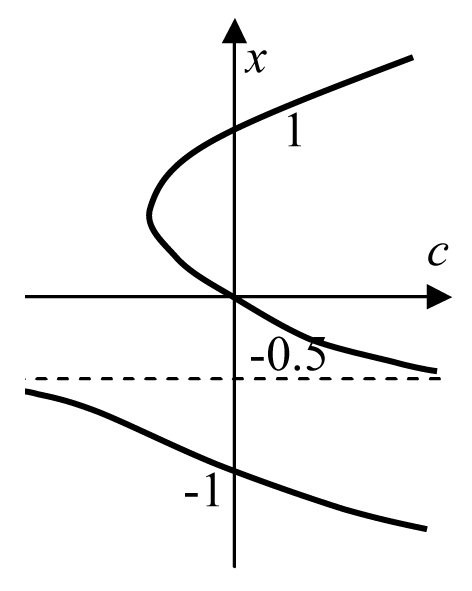}\\
a&b&c
\end{array}$\\
$\begin{array}{c}
\includegraphics[width=2.7cm]{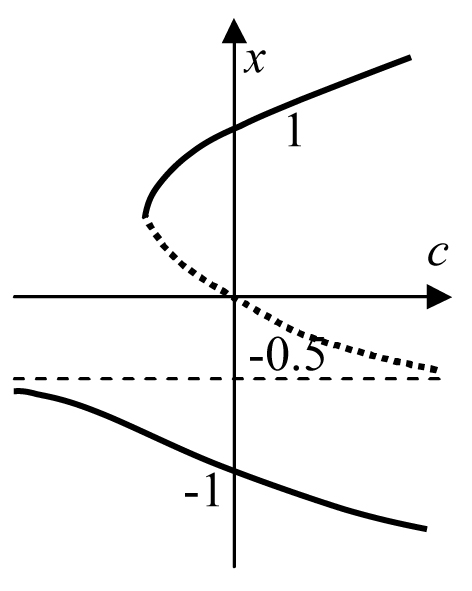}\\
d
\end{array}$
\caption{Example 2-a}\label{Ex2a}
\end{center}
\end{figure}

\begin{figure}[h]
\begin{center}$
\begin{array}{cc}
\includegraphics[width=2.7cm]{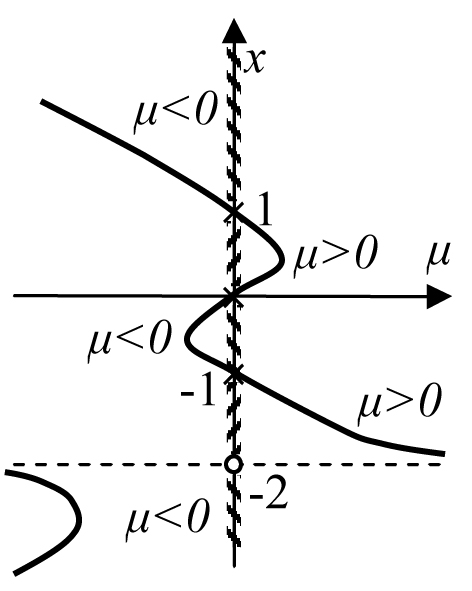}&
\includegraphics[width=2.7cm]{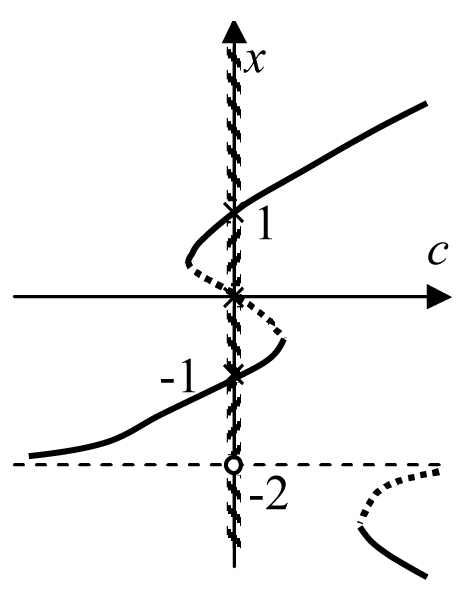}\\
a&b
\end{array}$
\caption{Example 2-b}\label{Ex2b}
\end{center}
\end{figure}

\begin{Example}
Consider dynamical system $\dot x=c+dx-x^3$.
\end{Example}
a) Plotting the bifurcation diagram for $d=1+2c$.
\newline
Step1: $c+(1+2c)x-x^3 \Rightarrow 1-c \frac{2x+1}{x^3-x}=1+\mu \frac{2x+1}{x(x^2-1)}$, where $\mu=-c$. There are three poles at $x=0, 1, -1$, and one zero at $x=-\frac{1}{2}$.
Step 2: See Fig.\ref{Ex2a}.a.
Step 3: There is no finite value vertical asymptote.
Step 4: Consider a horizontal asymptote at $x=-0.5$.
Step 5: See Fig.\ref{Ex2a}.b.
Step 7:    See Fig.\ref{Ex2a}.c.
Step 8: See Fig.\ref{Ex2a}.d.
\newline
b) Plotting the bifurcation diagram for $d=1+0.5c$. See Fig.\ref{Ex2b}.

\begin{figure}[h]
\begin{center}$
\begin{array}{ccc}
\includegraphics[width=2.5cm]{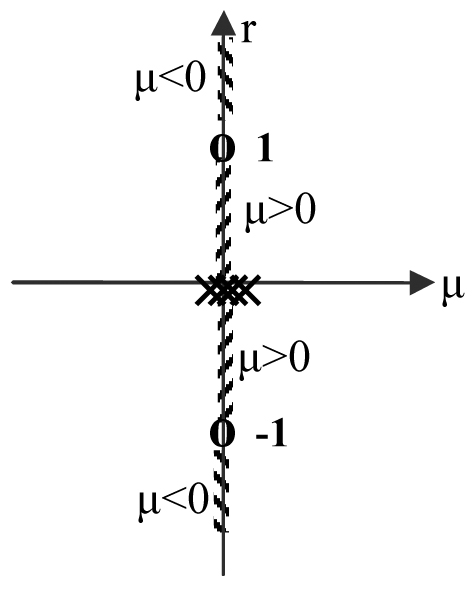}&
\includegraphics[width=2.5cm]{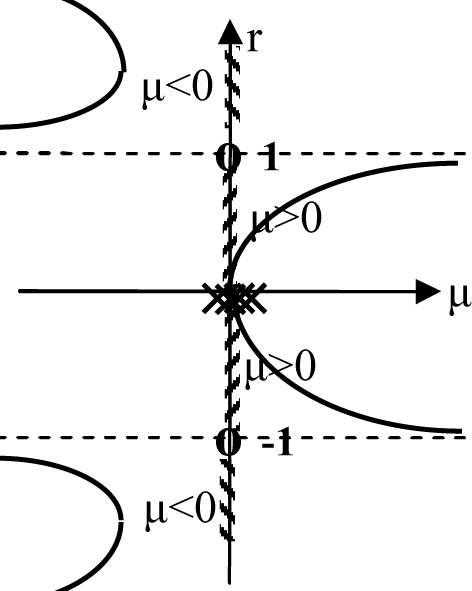}&
\includegraphics[width=2.5cm]{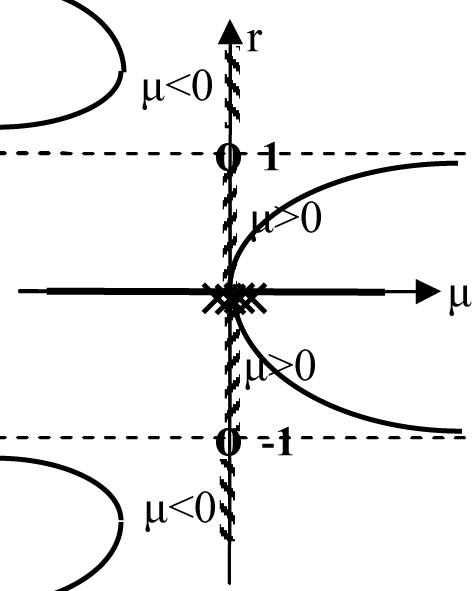}\\
a&b&c
\end{array}$ \\
$\begin{array}{cc}
\includegraphics[width=2.7cm]{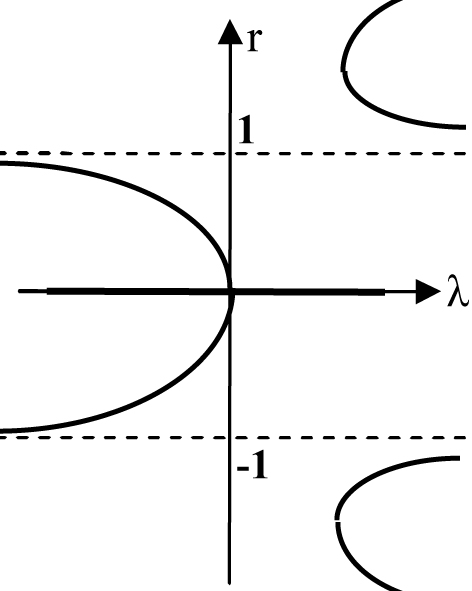}&
\includegraphics[width=2.7cm]{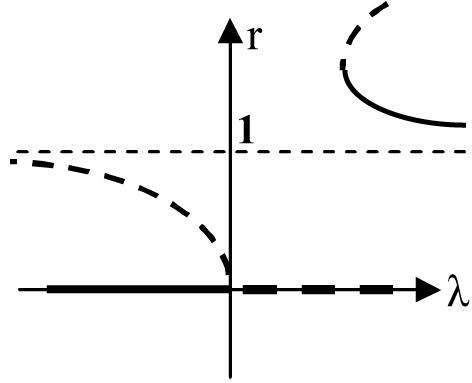}\\
d&e
\end{array}$
\caption{Example 3}\label{Ex3}
\end{center}
\end{figure}

\begin{Example}
Consider dynamical system
\begin{eqnarray}
\dot x &=& y+x[\lambda-\eta(x^2+y^2)+(x^2+y^2)^2] \nonumber\\
\dot y &=& -x+y[\lambda-\eta(x^2+y^2)+(x^2+y^2)^2] \nonumber
\end{eqnarray}
Plotting the bifurcation diagram for $\lambda=\eta$.
\end{Example}
Rewrite the system in polar coordination with $\lambda=\eta$:
\begin{eqnarray}
\dot r &=& r[\lambda-\lambda r^2+r^4] \nonumber\\
\dot \theta &=& -1 \nonumber
\end{eqnarray}
Step1: $\dot r=r[1-\lambda \frac{r^2-1}{r^4}]=r[1+\mu \frac{r^2-1}{r^4}]$, where $\mu=-c$.
Step 2: See Fig.\ref{Ex3}.a.
Step 3: There is no finite value vertical asymptote.
Step 4: Consider two horizontal asymptotes at $x=-1, 1$.
Step 5: See Fig.\ref{Ex3}.b.
Step 6: See Fig.\ref{Ex3}.c.
Step 7:    See Fig.\ref{Ex3}.d.
Step 8: Since $r>0$, only the upper half plane must be considered. See Fig.\ref{Ex3}.e.

\section{Conclusion}
\label{Conclusion}
The proposed algorithm prepares a simple and fast method to plot the local bifurcation diagram. This algorithm can be applied to dynamical systems introduced by Eq.1. It consists of some handy steps. Some examples show its performance and properties in the sense of simplification and quickness.


%

\end{document}